\documentclass{amsart}
\usepackage{enumitem,amssymb}
\usepackage[margin=1in]{geometry}
\usepackage[all]{xy}
\usepackage{xcolor}
\usepackage{mathtools}
\usepackage{epsdice}

\usepackage{relsize}
\usepackage[bbgreekl]{mathbbol}
\usepackage{amsfonts}
\DeclareSymbolFontAlphabet{\mathbb}{AMSb} 
\DeclareSymbolFontAlphabet{\mathbbl}{bbold}
\newcommand{\Prism}{{\mathlarger{\mathbbl{\Delta}}}}

\usepackage{enumitem}
\usepackage[colorlinks=true,hyperindex, linkcolor=magenta, pagebackref=false, citecolor=cyan,pdfpagelabels]{hyperref}
\usepackage[capitalize]{cleveref}
\setlist[enumerate]{itemsep=2pt,parsep=2pt,before={\parskip=2pt}}

\newtheorem{theorem}{Theorem}[section]
\newtheorem{proposition}[theorem]{Proposition}

\theoremstyle{definition}

\newtheorem{remark}[theorem]{Remark}

\newtheorem{notation}[theorem]{Notation}

\newtheorem*{convention*}{Conventions}

\begin{document}

\title{Totaro's inequality for classifying spaces}
\begin{abstract}
For a complex Lie group $G$ and a prime number $p$, Totaro had conjectured that the dimension of the singular cohomology with $\mathbf{Z}/p$-coefficients of classifying space of $G$ is bounded above by that of the de Rham cohomology of the classifying stack of (the split form of) $G$ in characteristic $p$.  This conjecture was recently proven by Kubrak--Prikhodko. In this note, we give a shorter proof.
\end{abstract}

\author{Bhargav Bhatt}
\author{Shizhang Li}
\maketitle

\section{The main theorem}

The goal of this note is to give a short proof of the following result \cite[Theorem 5.6.3]{KPHodge}, giving an algebro-geometric upper bound for the mod $p$ Betti numbers of the classifying space $BG(\mathbf{C})$:

\begin{theorem}[Kubrak--Prikhodko]
\label{mainthm}
Fix a split reductive group $G/\mathbf{Z}$. For any prime number $p$, we have an inequality
\begin{equation}
\label{eq:TotIneq} \dim H^i(BG_{\mathbf{C}}; \mathbf{Z}/p) \leq \dim H^i_{dR}(BG_{\mathbf{F}_p})
\end{equation}
for all $i$. 

\end{theorem}

The inequality \eqref{eq:TotIneq} continues a line of inquiry initiated by Totaro in \cite{TotaroHodgeClassifying}. In that paper, Totaro proved that \eqref{eq:TotIneq} was an equality for ``non-torsion'' primes of $G$, gave examples to show why it cannot be an equality in general, and conjectured the general case. Following further special cases exhibited by Primozic \cite{PrimozicHodgeClassifying}, the inequality was proven in full generality \cite[Theorem 5.6.3]{KPHodge}. 

The proof of \eqref{eq:TotIneq} in \cite{KPHodge} involves two steps. First, as anticipated in \cite{TotaroHodgeClassifying}, they extend the prismatic techniques \cite{BMS1,BhattScholzePrism} (which give analogues of inequality \eqref{eq:TotIneq} for smooth projective varieties) to stacks \cite[\S 2]{KPHodge}; this reduces inequality \eqref{eq:TotIneq}  to a purely rigid analytic statement over $\mathbf{Q}_p$ involving a comparison of the $\mathbf{F}_p$-cohomology of two classifying stacks that one can build starting from a reductive group scheme $G/\mathbf{Z}_p$ (corresponding to the algebraic group $G_{\mathbf{C}_p}$ and its ``maximal compact'' subgroup $\widehat{G}_{\mathbf{C}_p}$). Next, they prove this comparison \cite[Theorem 5.6.2]{KPHodge} by using the structure theory of reductive algebraic groups to devissage down to $G=\mathbf{G}_m$, where one can check the statement explicitly.  In this note, we observe that the second step of \cite{KPHodge} can be simplified (at the expense of proving a slightly weaker statement that nevertheless suffices for Theorem~\ref{mainthm}, see Remark~\ref{KPeqrmk}) by approximating the classifying stack $BG_{\mathbf{C}_p}$ by a smooth projective variety with good reduction at $p$ (Proposition~\ref{KPineq}).

\subsection*{Acknowledgements}  We thank Dmitry Kubrak for conversations and comments, and Monahan's for providing inspiring work conditions. B.B. was was partially supported by the NSF (\#1801689,\ \#1952399,\ \#1840234), a Packard fellowship, and the Simons Foundation (\#622511).

\section{Proof}

\begin{notation}
Fix a prime $p$. Let $\mathcal{O}_K$ be a ring of integers in a $p$-adic field $K$ with algebraically closed residue field $k$. Write $C = \widehat{\overline{K}}$ for a completed algebraic closure of $K$; this is an algebraically closed non-archimedean field with valuation ring $\mathcal{O}_C$. Write $C^\flat := \varprojlim_{x \mapsto x^p} C$ for the tilt of $C$ in the sense of Fontaine; this multiplicative monoid is naturally an algebraically closed field of characteristic $p$. 

Write $X \mapsto X^{an}$ for the analytification functor on varieties over $K$ (resp.~$C$), viewed as landing in adic spaces over $\mathrm{Spa}(K, \mathcal{O}_K)$ (resp.~$\mathrm{Spa}(C,\mathcal{O}_C)$). Likewise, write $\mathfrak{X} \mapsto \mathfrak{X}_\eta$ for the functor carrying a formal scheme $\mathfrak{X}$ over $\mathrm{Spf}(\mathcal{O}_K)$ (resp.~$\mathrm{Spf}(\mathcal{O}_C)$) to adic generic fibre $\mathfrak{X}_\eta$ over $\mathrm{Spa}(K, \mathcal{O}_K)$ (resp.~$\mathrm{Spa}(C,\mathcal{O}_C)$). We shall often use the same notation for the analogous functors for algebraic stacks over $K$ (resp.~$C$) or formal algebraic stacks over $\mathcal{O}_K$ (resp.~$\mathcal{O}_C$) (for the \'etale topology). In this case, to avoid any subtleties, we merely regard the output as a sheaf of groupoids for the pro-\'etale topology on classical rigid spaces over $K$ (resp.~$C$); this suffices for our purposes as we merely need access to the \'etale cohomology of these analytifications, and nothing more refined.

We shall implicitly use certain comparison theorems in \'etale cohomology, so let us recall the statements. For a finite type $C$-scheme $X$, we can identify the usual \'etale cohomology groups $H^*(X;\mathbf{Z}/p)$ of $X$ naturally with the \'etale cohomology groups $H^*(X^{an};\mathbf{Z}/p)$ of the adic space $X^{an}$, thanks to \cite[Theorem 3.8.1]{HuberEtaleBook}. Moreover, if one fixes an isomorphism $C \simeq \mathbf{C}$, then we can also identify these groups with the singular cohomology groups $H^*(X(\mathbf{C});\mathbf{Z}/p)$ of the complex analytic space $X(\mathbf{C})$ by the Artin comparison theorem. Both these isomorphisms extend to finite type algebraic stacks over $C$ by descent. 

\end{notation}

\begin{proof}[Proof of~Theorem~\ref{mainthm}]
Fix a reductive group scheme $G/\mathcal{O}_K$. We shall show 
\[ \dim_{\mathbf{F}_p} H^i(BG_{{C}}; \mathbf{Z}/p) \leq \dim H^i_{dR}(BG_{k}),\]
which suffices by standard comparison theorems. Using Proposition~\ref{KPineq} below, it is enough to show
\[ \dim_{\mathbf{F}_p} H^i(\widehat{BG}_{\eta,C}; \mathbf{Z}/p) \leq  \dim_{k} H^i_{dR}(BG_k). \]
For notational ease, write $\mathfrak{X} := \widehat{BG}$ for the formal completion of $BG$. We shall deduce the above inequality from the formal properties of the the prismatic cohomology of $\mathfrak{X}/\mathcal{O}_K$ (relative to the so-called Breuil--Kisin prism), mimicking the argument for smooth proper formal $\mathcal{O}_K$-schemes in \cite{BMS1,BhattScholzePrism}.

First, we introduce some prismatic notation. Choose a uniformizer $\pi \in \mathcal{O}_K$, yielding a surjection $\mathfrak{S} \coloneqq W(k)[\![u]\!] \twoheadrightarrow \mathcal{O}_K$ with kernel $I$. Endowing $\mathfrak{S}$ with the Frobenius lift sending $u$ to $u^p$ and acting via the (unique) Frobenius lift on $W \subset \mathfrak{S}$, we obtain a prism $(\mathfrak{S},I)$ over $\mathcal{O}_K$.
Finally, choosing a compatible system $\{\pi^{1/p^n}\}$ of $p$-power roots of $\pi$ in $\mathcal{O}_C$, we get a natural map $ \mathfrak{S}/p = k[\![u]\!] \to C^{\flat}$ with $u \mapsto (\pi^{1/p^n})_{n}$; in fact, this map realizes $C^\flat$ as a completed algebraic closure of $k(\!(u)\!)$.

Using the prism $(\mathfrak{S},I)$, we define the prismatic cohomology of $\mathfrak{X}$ via descent as follows. Let $U^\bullet$ be the standard simplicial smooth affine scheme $U^\bullet$ presenting $BG$, so $U^n = G^n$; set
\[ R\Gamma_\Prism(\mathfrak{X}/\mathfrak{S}) := R\lim R\Gamma_\Prism(\widehat{U^\bullet}/\mathfrak{S}) \in D(\mathfrak{S}),\]
where the derived limit is over the simplex category (i.e., it is a totalization of the corresponding cosimplicial object). This object is independent of the choice of the smooth presentation, but we do not need this here. The resulting (mod $p$) prismatic cohomology complex
\[ M \coloneqq R\Gamma_{\Prism}(\mathfrak{X}/\mathfrak{S}) \otimes^L_{\mathfrak{S}} k[\![u]\!] \in D(k[\![u]\!])\] 
has the following properties:
\begin{enumerate}
\item $M$ is a derived $u$-complete object in $D^{\geq 0}(k[\![u]\!])$, with a semilinear Frobenius map
$\varphi \colon M \to M$.
\item One has a canonical isomorphism 
\[ (\varphi^* M) \otimes^L_{k[\![u]\!]} k  \simeq R\Gamma_{dR}(\mathfrak{X}_{k}) \quad \text{or equivalently} \quad M \otimes^L_{k[\![u]\!]} k  \simeq \varphi_* R\Gamma_{dR}(\mathfrak{X}_{k})\] 
in $D(k)$ relating $M$ with the de Rham cohomology of $\mathfrak{X}_{k}$ ($\simeq BG_{k}$).
\item After inverting $u$, the Frobenius map induces an isomorphism
$\varphi^*M[1/u] \xrightarrow{\simeq} M[1/u]$.
\item There is a canonical isomorphism 
\[R\Gamma(\mathfrak{X}_{\eta,C}; \mathbf{Z}/p) \simeq \left(M \otimes_{k[\![u]\!]} C^{\flat}\right)^{\varphi=1} \in D(\mathbf{Z}/p).\]
\end{enumerate}
These can be deduced via descent from similar properties for the smooth affine formal schemes comprising $U^\bullet$ that were proven in \cite[Theorem 1.8.(4)-(6)]{BhattScholzePrism}. Indeed, (1) follows as derived $u$-completeness is preserved under limits; (2) follows from smooth descent for de Rham cohomology (see \cite[\S 2]{ABM}) and implies (via derived Nakayma) that each cohomology group $H^i(M)$ is a finitely generated $k[\![u]\!]$-module by the finiteness for de Rham cohomology in each degree shown in \cite{TotaroHodgeClassifying} (argument recalled in next paragraph); using this finiteness to avoid extra completions, (3) and (4) then follow as totalizations in $D^{\geq 0}$ commute with filtered colimits and flat base change (such as $\varphi^*$, inverting $u$, and tensoring up to $C^\flat$) and because $R\Gamma(\mathfrak{X}_{\eta,C};\mathbf{Z}/p) \simeq R\lim R\Gamma(\widehat{U^\bullet}_{\eta,C}; \mathbf{Z}/p)$ by descent.

We now prove the desired inequality. Pick an arbitrary integer $i$.
Property (1) implies that each $H^i(M)$ is a derived $u$-complete module over $k[\![u]\!]$.
Property (2) implies that we have an injection 
\[
\label{specialization}
\tag{\epsdice{1}}
H^i(M)/u \hookrightarrow H^i(M/u) = \varphi_* H^i_{dR}(BG_k),
\]
with the target being a finite dimensional $k$-vector space due to~\cite{TotaroHodgeClassifying} (and the perfectness of $k$).
These two altogether imply that each $H^i(M)$ is a finitely generated module over $k[\![u]\!]$: any map $k[\![u]\!]^{\oplus n} \to H^i(M)$ that is surjective modulo $u$ must be surjective as the cokernel, being both derived $u$-complete and $u$-divisible, is zero. The finiteness plus property (3) imply that (c.f.~\cite[Lemma 8.5]{B15}) the natural map
\[\left(H^i(M) \otimes_{k[\![u]\!]} C^{\flat}\right)^{\varphi=1} \otimes_{\mathbf{F}_p} C^{\flat} \rightarrow H^i(M) \otimes_{k[\![u]\!]} C^{\flat}\]
is an isomorphism.
In particular, together with property (4), we have
\[
\dim_{\mathbf{F}_p} H^i(\mathfrak{X}_{\eta,C}; \mathbf{Z}/p) = \dim_{C^{\flat}} H^i(M) \otimes_{k[\![u]\!]} C^{\flat}
= \dim_{k(\!(u)\!)} H^i(M)[1/u].
\]
Since $H^i(M)$ is finitely generated over $k[\![u]\!]$, using~\eqref{specialization}, we have inequalities
\[
\dim_{k(\!(u)\!)} H^i(M)[1/u] \leq \dim_{k} H^i(M)/u \leq \dim_{k} H^i(M/u) = \dim_{k} H^i_{dR}(BG_k).
\]
Combining these gives the desired inequality
\[ \dim_{\mathbf{F}_p} H^i(\mathfrak{X}_{\eta,C}; \mathbf{Z}/p) \leq  \dim_{k} H^i_{dR}(BG_k). \qedhere \]

\end{proof}

We have reduced our main theorem to the following assertion:

\begin{proposition}
\label{KPineq}
Let $k$ be an algebraically closed field of characteristic $p$, and let $W = W(k)$. Fix a split reductive group scheme $G/W$. Then 
\[ \dim H^i(BG_C; \mathbf{Z}/p) \leq \dim H^i(\widehat{BG}_{\eta,C}; \mathbf{Z}/p)\]
for all $i$, where $C$ denotes a completed algebraic closure of $W[1/p]$. 
\end{proposition}

\begin{remark}
\label{KPeqrmk}
Kubrak--Prikhodko have shown that the inequality in Proposition~\ref{KPineq} is in fact always an equality (see \cite[Theorem 5.6.2]{KPHodge}, which applies to a more general situation). Their result can be regarded roughly as a $p$-adic analog of the fact that the homotopy type does not change when one passes from a complex Lie group to its maximal compact subgroup.
\end{remark}

We prove Proposition~\ref{KPineq} by approximating the cohomology of $BG_C$ carefully. Specifically, Ekedahl \cite{EkedahlApprox} has constructed a smooth projective variety $Z/C$ and a map $Z \to BG_C$ that induces an injection on cohomology in a given range of degrees.  By constructing $Z$ with some care given to integral models,  we show that the analytification of Ekedahl's map can be factored over the canonical map $\widehat{BG}_{\eta,C} \to BG_C^{an}$, which immediately yields the inequality.  Implementing this strategy requires some basic notions of GIT over $W$; we refer to Seshadri's paper \cite{SeshadriGeomRed} for background.

\begin{proof}
Fix an integer $N \geq 0$. We shall show the desired inequality for all $i \leq N$. Letting $N \to \infty$ then implies the proposition. 

Let $H = G \times \mathbf{G}_m$.  Choose a faithful representation $H \to \mathrm{GL}(V)$ over $W$ which restricts to scalar multiplication on the summand $\mathbf{G}_m$ of $H$ (so $G$ acts on $\mathbf{P}(V)$). By replacing $V$ with $V^{\oplus n}$ for $n \gg 0$, we may assume that there is an $H$-stable open $U \subset V$ whose complement has arbitrarily high codimension $c'$ (relative to our choice of $N$), and such that $H$ acts freely on $U$. Killing the scalars shows that $U' := U/\mathbf{G}_m \subset \mathbf{P}(V)$ carries a free $G$-action (as $G = H/\mathbf{G}_m$), lies in the $G$-semistable locus $\mathbf{P}(V)^{ss}$ (in fact, it lives in the stable locus), and the closed set $\mathbf{P}(V)^{ss} - U' \subset \mathbf{P}(V)^{ss}$ has codimension $\geq c'$. 

Let $X := \mathbf{P}(V)//G = \mathrm{Proj}(\mathrm{Sym}^*(V)^G)$ be the GIT quotient. Standard results in GIT over $W$ (see \cite[Theorem 4]{SeshadriGeomRed}) show that there is a surjective $G$-invariant quotient map $q:\mathbf{P}(V)^{ss} \to X$, and that $X$ is a normal projective $W$-scheme (in fact, the $k$-th power of $\mathcal{O}_{\mathbf{P}(V)^{ss}}(1)$ descends along $q$ for some $k \geq 1$). The $G$-stable open $U' \subset \mathbf{P}(V)^{ss}$ descends to a smooth open $\overline{U} \subset X$, and the quotient map $q:U' \to \overline{U}$ is a $G$-torsor. The complement $X - \overline{U} \subset X$ is the image of $\mathbf{P}(V)^{ss} - U' \subset \mathbf{P}(V)^{ss}$, and hence has codimension $\geq c \coloneqq c' - \dim G$ in $X$. 

Let $Q/W$ be the projective space parametrizing complete intersections of dimension $c-1$ in $X$. Ekedahl has shown in \cite[\S 1.1]{EkedahlApprox} that, for $c \gg 0$, there exists a dense Zariski open subset $Q_{1} \subset Q_{W[1/p]}$ such that  any complete intersection $Z \subset X_{C}$ corresponding to a point of $Q_{1}(C)$ is smooth, lies in $\overline{U}_C \subset X_{C}$, and has the property that the map $H^i(BG_C;\mathbf{Z}/p) \to H^i(Z; \mathbf{Z}/p)$ induced by restricting the $G_C$-torsor $q$ is injective for $i \leq N$. On the other hand, consider the fibrewise dense Zariski open $Q_2 \subset Q$ parametrizing smooth complete intersections that lie in $\overline{U}$. We first observe that $Q_2(\mathcal{O}_C) \cap Q_{1}(C) \subset Q(C)$ is non-empty for general reasons: a proper Zariski closed subset (such as the analytification of $Q_{C} - Q_{1,C}$) of the smooth proper connected rigid space $Q_{C}^{an}$ cannot contain non-empty quasi-compact open subspaces (such as $(\widehat{Q_{2,\mathcal{O}_C}})_\eta \subset Q_{C}^{an}$). A point of $Q_2(\mathcal{O}_C) \cap Q_{1}(C) \subset Q(C)$ then gives a smooth complete intersection $Z \subset X_{\mathcal{O}_C}$ that actually lies in $\overline{U}_{\mathcal{O}_C}$ and such that the map $H^i(BG_C;\mathbf{Z}/p) \to H^i(Z_C; \mathbf{Z}/p)$ considered above is injective for $i \leq N$. Fix one such smooth complete intersection $Z$. As $Z \subset \overline{U}_{\mathcal{O}_C}$, the map $q$ restricts to a $G_{\mathcal{O}_C}$-torsor over $Z$ (and not merely $Z_C$). By functoriality passing to the adic generic fibre before and after $p$-adic completions, this gives a commutative diagram
\[ \xymatrix{  \widehat{Z}_\eta \ar[r] \ar[d] & Z_C^{an} \ar[r] \ar[d] & Z_C \ar[d] \\
			\widehat{BG}_{\eta,C} \simeq (\widehat{BG_{\mathcal{O}_C}})_\eta \ar[r] & (BG_C)^{an} \ar[r] & BG_C.}\]
where the right vertical map comes as the scheme-theoretic generic fibre of $Z \to BG_{\mathcal{O}_C}$, the middle vertical map is the analytification of the right vertical map, and the left vertical map is the adic generic fibre of the $p$-adic completion of $Z \to BG_{\mathcal{O}_C}$. Now both of the top horizontal maps induce an isomorphism on $H^i(-;\mathbf{Z}/p)$: the left one is an isomorphism by properness, while the right one gives an isomorphism on cohomology by Huber's comparison theorem. As the right vertical map is injective on $H^i(-;\mathbf{Z}/p)$ for $i \leq N$ by construction, it follows that the same holds for the composite bottom horizontal maps, which gives the desired claim.
\end{proof}



\begin{thebibliography}{BS}

\bibitem[ABM]{ABM}
B.~Antieau, B.~Bhatt, A.~Mathew, {\em Counterexamples to Hochschild-Kostant-Rosenberg in characteristic $p$}, arXiv:1909.11437, to appear in {\em Forum of Mathematics, Sigma}.

\bibitem[B]{B15}
B.~Bhatt, {\em Specializing varieties and their cohomology from characteristic 0 to characteristic {$p$}}, Algebraic geometry: Salt Lake City 2015, 43--88,
Proc.~Sympos.~Pure Math., 97.2, Amer.~Math.~Soc., Providence, RI, 2018.

\bibitem[BMS]{BMS1}
B.~Bhatt, M.~Morrow, P.~Scholze, {\em Integral $p$-adic Hodge theory}, Publ.~Math.~Inst.~Hautes \'{E}tudes Sci.~128 (2018), 219--397. 

\bibitem[BS]{BhattScholzePrism}
B.~Bhatt, P.~Scholze, {\em Prisms and prismatic cohomology}, arXiv:1905.08229

\bibitem[E]{EkedahlApprox}
T.~Ekedahl, {\em Approximating classifying spaces by smooth projective varieties}, arXiv: 0905.1538

\bibitem[H]{HuberEtaleBook}
R.~Huber, {\em \'Etale cohomology of rigid-analytic varieties and spaces}, Aspects of Mathematics, E30. Friedr. Vieweg \& Sohn, Braunschweig, 1996

\bibitem[KP]{KPHodge}
D.~Kubrak and A.~Prikhodko, {\em $p$-adic Hodge theory for stacks}, arXiv: 2105.05319

\bibitem[P]{PrimozicHodgeClassifying}
E.~Primozic, {\em Computations of de Rham cohomology rings of classifying stacks at torsion primes}, arXiv: 1909.13413

\bibitem[S]{SeshadriGeomRed}
C.~Seshadri, {\em Geometric reductivity over arbitrary bases},  Advances in Math.~26 (1977), no.~3, 225--274. 

\bibitem[T]{TotaroHodgeClassifying}
B.~Totaro, {\em Hodge theory for classifying stacks}, Duke Math.~J.~167 (2018), no.~8, 1573--1621. 



\end{thebibliography}
\end{document}